\documentclass{ifacconf}

\usepackage{graphicx}      
\usepackage{natbib}        
\begin{document}
\begin{frontmatter}

\title{Continuation model predictive control on smooth manifolds} 

\thanks[footnoteinfo]{Accepted to the 16th IFAC Workshop on Control Applications
of Optimization (CAO'2015), Garmisch-Partenkirchen, Germany, October 6--9, 2015.}

\author[First]{Andrew Knyazev} 
\author[Second]{Alexander Malyshev}

\address[First]{Mitsubishi Electric Research Labs (MERL)
201 Broadway, 8th floor, Cambridge, MA 02139, USA. (e-mail: knyazev@merl.com), \hfill
({http://www.merl.com/people/knyazev}).}
\address[Second]{Mitsubishi Electric Research Labs (MERL)
201 Broadway, 8th floor, Cambridge, MA 02139, USA. (e-mail: malyshev@merl.com)}
\begin{abstract}                
Model predictive control (MPC) anticipates future events to take appropriate control actions.
Nonlinear MPC (NMPC) describes systems with nonlinear models and/or constraints.
Continuation MPC, suggested by T.~Ohtsuka in 2004, uses Krylov-Newton iterations.
Continuation MPC is suitable for nonlinear problems and has been recently adopted
for minimum time problems. We extend the continuation MPC approach to a case
where the state is implicitly constrained to a smooth manifold.
We propose an algorithm for on-line controller implementation and demonstrate
its numerical effectiveness for a test problem on a hemisphere.
\end{abstract}


\end{frontmatter}

\section{Introduction}

Model predictive control (MPC) is exposed, e.g., in \cite{CaBo:04,GrPa:11}. MPC is
efficient in industrial applications; see \cite{QiBa:03}. Numerical aspects of MPC are
discussed in \cite{DiFeHa:09}. Our contribution to numerical methods for MPC is
further development of the approach by \cite{Oht:04} and efficient use of
the Newton-Krylov approximation for numerical solving the nonlinear MPC problems.
\cite{KnFuMa:15} solves a minimum-time problem and investigates preconditioning
for the Newton-Krylov method.

The MPC method owes its success to efficient treatment of constraints on control
and state variables. The goal of the present note is to draw attention to the important
fact that, in addition to the efficient treatment of explicit constraints, MPC can easily
incorporate the structure-preserving geometric integration of ordinary differential equations
modeling the system.

The idea of structure-preserving numerical methods for differential equations has been employed
in numerical analysis since 1950s, when the first computers began to be used in engineering
computations. However, its systematic study under the name of
\emph{geometric integration of differential equations} has been accomplished
within the last 40 years. An accessible introduction to this active research domain is found in \cite{Hai:11}. A more advanced source is \cite{HaWaLu:06}.

Variational formulations of many models described by differential equations
lead to the Hamiltonian system dynamics; see \cite{LeRe:04}. Such dynamic systems
conserve the energy and possess a special geometric structure called the symplectic structure.
Much effort has been spent by numerical analysts to develop the called symplectic
numerical methods, which are especially efficient for long-time integration,
because they produce qualitatively correct computed solutions.
The symplectic methods are superior even in the short-time simulations.
\cite{LeRe:04} presents the state-of-the-art symplectic numerical algorithms
and supporting theory.

Description of the system dynamics in the form of ordinary differential
equations appears twice in MPC: first, within the finite-horizon prediction problem and,
secondly, when advancing the current state with the input control computed by
the finite-horizon prediction. The latter can be omitted if the states of the system
are measured directly by the controller sensors. Furthermore, the application of the geometric
structure-preserving integration algorithms from \cite{HaWaLu:06} or other sources for advancing
the current state is straightforward. 

The use of the geometric integration during the finite-horizon prediction is less obvious.
The main difference between the method from \cite{KnFuMa:15} and the variant of the present paper
is in that the latter applies a structure-preserving geometric integration
inside the forward recursion.

The rest of the note is organized as follows. Section 2 presents a framework of
the finite-horizon prediction problem and expresses its solution in the form of
a nonlinear algebraic equation. Geometric integration is incorporated during
the elimination of state variables from the KKT conditions. Section 3 discusses how classical
integration methods can be used for integration on manifolds. The content of Section 3
is mainly adopted from \cite{Hai:01} and discusses the local coordinates approach and
projection methods. Section 4 describes a simple test example and necessary formulas for
computer implementation. Section 5 shows numerical results and plots.

\section{Finite-horizon prediction}

Our model finite-horizon control problem, see \cite{KnFuMa:15}, along a fictitious time
$\tau\in[t,t+T]$ consists in choosing the control $u(\tau)$ and parameter vector $p$,
which minimize the performance index $J$ as follows:
\[
\min_{u,p} J,
\]
where
\[
J = \phi(x(t+T),p)+\int_t^{t+T}L(\tau,x(\tau),u(\tau),p)d\tau
\]
subject to the equation for the state dynamics
\begin{equation}\label{e1}
\frac{dx}{d\tau}=f(\tau,x(\tau),u(\tau),p),
\end{equation}
and the equality constraints for the state $x$ and the control~$u$
\begin{equation}\label{e2}
C(\tau,x(\tau),u(\tau),p) = 0,
\end{equation}
\begin{equation}\label{e3}
\psi(x(t+T),p) = 0.
\end{equation}
The initial value condition $x(\tau)|_{\tau=t}$ for (\ref{e1}) is the state vector $x(t)$
of the dynamic system. The control vector $u=u(\tau)|_{\tau=t}$ is used afterwards as an input
to control the system at time~$t$.
The components of the vector $p(t)$ are parameters of the dynamic system.
The horizon time length $T$ may depend on $t$.

The continuous formulation of the finite-horizon prediction problem stated above is discretized
on a time grid $\tau_i$, $i=0,1,\ldots,N$, through the horizon $[t,t+T]$ partitioned into
$N$ time steps of size $\Delta\tau_i=\tau_{i+1}-\tau_i$, and the time-continuous vector
functions $x(\tau)$ and $u(\tau)$ are replaced by their indexed values $x_i$ and $u_i$
at the grid points $\tau_i$. The integral of the performance cost $J$ over the horizon
is approximated by the rectangular quadrature rule. Equation (\ref{e1}) is approximated
by a one-step integration formula, for example, a Runge-Kutta method; cf. \cite{Hai:01}.
The discretized optimal control problem is as follows:
\[
\min_{u_i,p}\left[
\phi(x_N,p) + \sum_{i=0}^{N-1}L(\tau_i,x_i,u_i,p)\Delta\tau_i\right],
\]
subject to
\begin{equation}\label{e4}
\quad x_{i+1} = x_{i}+\Phi_i(\tau_i,x_i,u_i,p)\Delta\tau_i,\quad i = 0,1,\ldots,N-1,
\end{equation}
\begin{equation}\label{e5}
C(\tau_i,x_i,u_i,p) = 0,\quad  i = 0,1,\ldots,N-1,
\end{equation}
\begin{equation}\label{e6}
\psi(x_N,p) = 0.
\end{equation}
The function $\Phi_i$ is implicitly determined by the structure-preserving
method used for numerical integration of (\ref{e1}). We note that $\Phi_i(\tau_i,x_i,u_i,p)=f(\tau_i,x_i,u_i,p)+O(\Delta\tau)$, where
$\Delta\tau=\max_i\Delta\tau_i$.

The necessary optimality conditions for the discretized finite horizon problem
are obtained by means of the discrete Lagrangian function
\begin{eqnarray*}
\mathcal{L}(X,U)&=&\phi(x_N,p)+\sum_{i=0}^{N-1}L(\tau_i,x_i,u_i,p)\Delta\tau_i\\
&&{}+\lambda_0^T[x(t)-x_0]\\&&{}+\sum_{i=0}^{N-1}\lambda_{i+1}^T
[x_i-x_{i+1}+\Phi_i(\tau_i,x_i,u_i,p)\Delta\tau_i]\\
&&{}+\sum_{i=0}^{N-1}\mu_i^TC(\tau_i,x_i,u_i,p)\Delta\tau_i+\nu^T\psi(x_N,p),
\end{eqnarray*}
where $X = [x_i\; \lambda_i]^T$, $i=0,1,\ldots,N$, and
$U = [u_i\; \mu_i\; \nu\; p]^T$, $i=0,1,\ldots,N-1$.
Here, $\lambda$ is the costate vector and $\mu$ is the Lagrange multiplier vector 
associated with constraint~(\ref{e5}). The terminal constraint (\ref{e6})
is relaxed by the aid of the Lagrange multiplier $\nu$. 

The Hamiltonian function is denoted by
\begin{eqnarray*}
\lefteqn{H(t,x,\lambda,u,\mu,p) = L(t,x,u,p)}\hspace*{8em}&& \\
&&{}+\lambda^T f(t,x,u,p)+\mu^T C(t,x,u,p).
\end{eqnarray*}

The necessary optimality conditions are the KKT stationarity conditions:
$\mathcal{L}_{\lambda_i}=0$, $\mathcal{L}_{x_i}=0$, $i=0,1,\ldots,N$,
$\mathcal{L}_{u_j}=0$, $\mathcal{L}_{\mu_j}=0$, $i=0,1,\ldots,N-1$,
$\mathcal{L}_{\nu_k}=0$, $\mathcal{L}_{p_l}=0$.

Since $\partial\Phi_i(\tau_i,x_i,u_i,p)/\partial x_i$ is not always available,
we suggest to use $\partial f(\tau_i,x_i,u_i,p)/\partial x_i$ instead.
This modification is applied in the backward recursion used below.

The KKT conditions are reformulated in terms of a mapping $F[U,x,t]$,
where the vector $U$ combines the control input $u$, the Lagrange multiplier
$\mu$, the Lagrange multiplier $\nu$, and the parameter $p$, all in one vector:
\[
U(t)=[u_0^T,\ldots,u_{N-1}^T,\mu_0^T,\ldots,\mu_{N-1}^T,\nu^T,p^T]^T. 
\]
The vector argument $x$ in $F[U,x,t]$ denotes the current measured or estimated 
state vector, which serves as the initial vector $x_0$ in the following procedure.
\begin{enumerate}
\item Starting from the current measured or estimated state $x_0$, compute
$x_i$, $i=0,1\ldots,N-1$, by the forward recursion
\[
x_{i+1} = x_i + \Phi_i(\tau_i,x_i,u_i,p)\Delta\tau_i.
\]
Then starting from
\[
\lambda_N=\frac{\partial\phi^T}{\partial x}(x_N,p)+
 \frac{\partial\psi^T}{\partial x}(x_N,p)\nu
\]
compute the costates $\lambda_i$, $i=N\!-\!1,\ldots,0$, by the backward recursion
\[
\lambda_i=\lambda_{i+1}+\frac{\partial H^T}{\partial x}
(\tau_i,x_i,\lambda_{i+1},u_i,\mu_i,p)\Delta\tau_i.
\]
\item Calculate $F[U,x,t]$, using just obtained $x_i$ and $\lambda_i$, as \pagebreak

\begin{eqnarray*}
\lefteqn{F[U,x,t]}\\
&&\hspace*{-2em}=\left[\begin{array}{c}\begin{array}{c}
\frac{\partial H^T}{\partial u}(\tau_0,x_0,\lambda_{1},u_0,\mu_0,p)\Delta\tau_0\\
\vdots\\\frac{\partial H^T}{\partial u}(\tau_i,x_i,\lambda_{i+1},u_i,\mu_i,p)\Delta\tau_i\\
\vdots\\\frac{\partial H^T}{\partial u}(\tau_{N-1},x_{N-1},\lambda_{N},u_{N-1},
\mu_{N-1},p)\Delta\tau_{N-1}\end{array}\\\;\\
\begin{array}{c}C(\tau_0,x_0,u_0,p)\Delta\tau_0\\
\vdots\\C(\tau_i,x_i,u_i,p)\Delta\tau_i\\\vdots\\
C(\tau_{N-1},x_{N-1},u_{N-1},p)\Delta\tau_{N-1}\end{array}\\\;\\
\psi(x_N,p)\\[2ex]
\begin{array}{c}\frac{\partial\phi^T}{\partial p}(x_N,p)+
\frac{\partial\psi^T}{\partial p}(x_N,p)\nu\\
+\sum_{i=0}^{N-1}\frac{\partial H^T}{\partial p}(\tau_i,x_i,
\lambda_{i+1},u_i,\mu_i,p)\Delta\tau_i\end{array}
\end{array}\right].
\end{eqnarray*}
\end{enumerate}
The equation with respect to the unknown vector $U(t)$
\begin{equation}\label{e7}
 F[U(t),x(t),t]=0
\end{equation}
gives the required necessary optimality conditions that are solved
on the controller in real time. All details about solving
(\ref{e7}) by the Newton-Krylov method are found in \cite{KnFuMa:15}.
We only recall that our method uses GMRES iterations for solving
linear systems with the Jacobian matrix $F_U$. To accelerate the
convergence of GMRES iterations, the matrix $F_U$ is computed exactly
at some time instances and then used as a preconditioner between
these time instances. 

\section{Geometric integration}

This section gives a short explanation of geometric integration used in our
illustrative example.
We consider the ordinary differential equations (\ref{e1}), which describe
the system dynamics. Suppose that a property $L(x(t))=const$ is fulfilled on
each solution $x(t)$ of (\ref{e1}), where the constant $const$ depends on the solution.
If the numerical method $x_{i+1} = x_{i}+\Phi_i(\tau_i,x_i,u_i,p)\Delta\tau_i$
satisfies the same property $L(x_i(\tau))=const$, we refer to it as
a structure preserving method. A notorious example is given by the spheres
$L(x(t))=\|x\|_2=const$.

The article by \cite{Hai:01} ``illustrates how classical integration methods
for differential equations on manifolds can be modified in order to preserve
certain geometric properties of the exact flow. Projection methods as well as
integrators based on local coordinates are considered.''

We adopt the simplest method from \cite{Hai:01}, which is called the
\emph{method of local coordinates} there. When it is feasible,
it is the most accurate of all structure preserving methods.
Let a manifold $\mathcal{M}$ contain the solution of $\dot{x}=f(x)$, which
we want to compute. Assume that $\alpha\colon U\to{R}^n$ is a local parametrization
of $\mathcal{M}$ near the state $x_i=\alpha(z_i)$. The coordinate change $x=\alpha(z)$
transforms the differential equation $\dot{x}=f(x)$ into
\begin{equation}\label{e8}
\alpha'(z)\dot{z}=f(\alpha(z)).
\end{equation}
This is an over-determined system of differential equations because the dimension of $z$
is less than $n=\dim x$. However, $f(x)$ is tangent to $\mathcal{M}$ by assumption,
therefore, (\ref{e8}) is equivalent to a system
\begin{equation}\label{e9}
\dot{z}=\beta(z), \quad z(\tau_i)=z_i.
\end{equation}
If the transfer from (\ref{e8}) to (\ref{e9}) is easy for implementation, then
the method of local coordinates is feasible.

The principal idea of the method of local coordinates is to perform one step of a numerical method
applied to (\ref{e9}) and to map the result via the transformation $\alpha$ back to the
manifold. One step $x_i\mapsto x_{i+1}$of the resulting algorithm is implemented as follows.

\textbf{Algorithm 1}{ (Local coordinates approach)}
\begin{itemize}
\item choose a local parametrization $\alpha$ and compute $z_i$ from $x_i=\alpha(z_i)$;
\item compute $\widehat{z}_{i+1}=z_i+\Phi_i(z_i)\Delta\tau_i$, the result of the one-step
method applied to (\ref{e9});
\item define the numerical solution by $x_{i+1}=\alpha(\widehat{z}_{i+1})$.
\end{itemize}

\subsection{Symmetric integration methods}

Mechanical systems often obey the Hamiltonian structure, and numerical
methods for such system must be symmetric or time-reversible in order to preserve
geometric properties of the Hamiltonian systems.
Assume that a one-step numerical method applied to a system $\dot{x}=f(t,x)$
over the interval $[t_i,t_{i+1}]$ produces the state $x_{i+1}$ from the state $x_i$.
The time-reversibility of the numerical method means that its application to the
time-reversed system $-\dot{x}=f(-t,x)$ over the interval $[-t_{i+1},-t_i]$
produces $x_i$ from $x_{i+1}$. This is a case for the trapezoidal rule, but not
for the explicit Euler method.

\cite{Hai:01} shows that symmetric methods perform qualitatively better for integration
over long time intervals. But most of the commonly used techniques for solving
differential equations on manifolds destroy the symmetry of underlying method.
To restore the symmetry, additional modifications of numerical methods are necessary.

We illustrate restoration of the symmetry for another standard technique of
geometric integration on manifolds called the \emph{projection methods}. For
an ordinary differential equation $\dot{y}=f(y)$, the one step integration $y_n\mapsto y_{n+1}$
proceeds as follows.

\textbf{Algorithm 2}{ (Standard projection method)}
\begin{itemize}
\item compute $\widehat{y}_{n+1}=\Phi_h(y_n)$ by any numerical integrator $\Phi_h$
applied to $\dot{y}=f(y)$, e.g. by a Runge-Kutta method;
\item project $\widehat{y}_{n+1}$ orthogonally onto the manifold $\mathcal{M}$
to obtain $y_{n+1}\in\mathcal{M}$.
\end{itemize}

\begin{figure}[ht]
\begin{center}
\includegraphics[width=8.4cm]{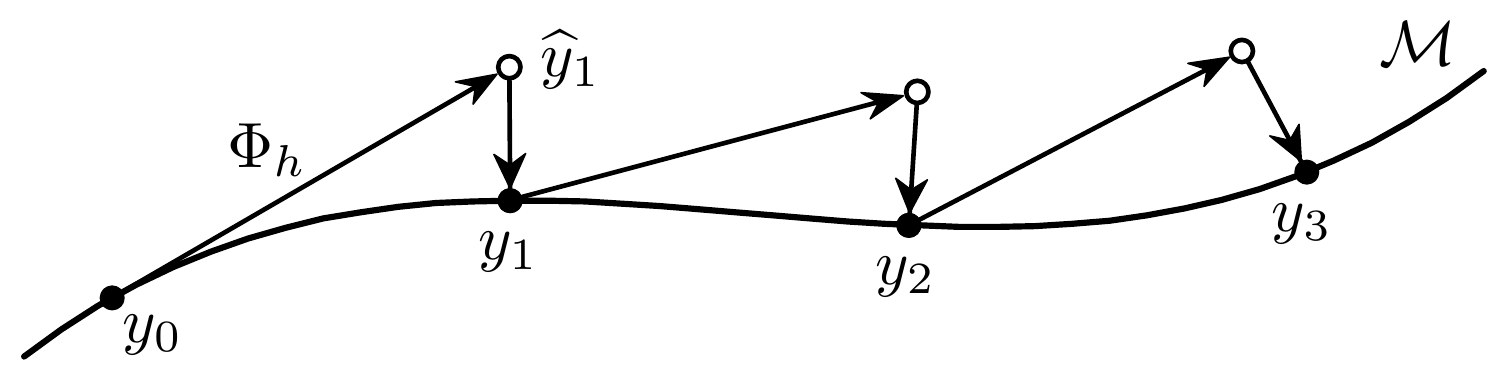}
\caption{The standard projection \cite[Fig. 3.1]{Hai:01}}
\label{figa}
\end{center}
\end{figure}

The standard projection method is illustrated in Figure~\ref{figa}.
The projection destroys the symmetry of $\Phi_h$ if it is available.
A symmetry-restoration algorithm is developed in \cite{Hai:00}.
The idea is to perturb the vector $y_n$ before applying a symmetric
one-step method such that the final projection is of the same size as the perturbation.

\textbf{Algorithm 3}{ (Symmetric projection method)}
\begin{itemize}
\item $\widetilde{y}_n=y_n+G^T(y_n)\mu$, where $g(y_n)=0$;
\item $\widehat{y}_{n+1}=\Phi_h(\widehat{y}_n)$ (symmetric one-step method for
 equation $\dot{y}=f(y)$);
\item $y_{n+1}=\widehat{y}_{n+1}+G^T(y_{n+1})\mu$ with vector $\mu$ such that
 $g(y_{n+1})=0$.
\end{itemize}
Here, $G(y)=g'(y)$ denotes the Jacobian of $g(y)$, if the manifold is given by
the condition $\mathcal{M}=\{y|g(y)=0\}$. It is important to choose the same vector $\mu$
in the perturbation an in the projection.

\begin{figure}[ht]
\begin{center}
\includegraphics[width=8.4cm]{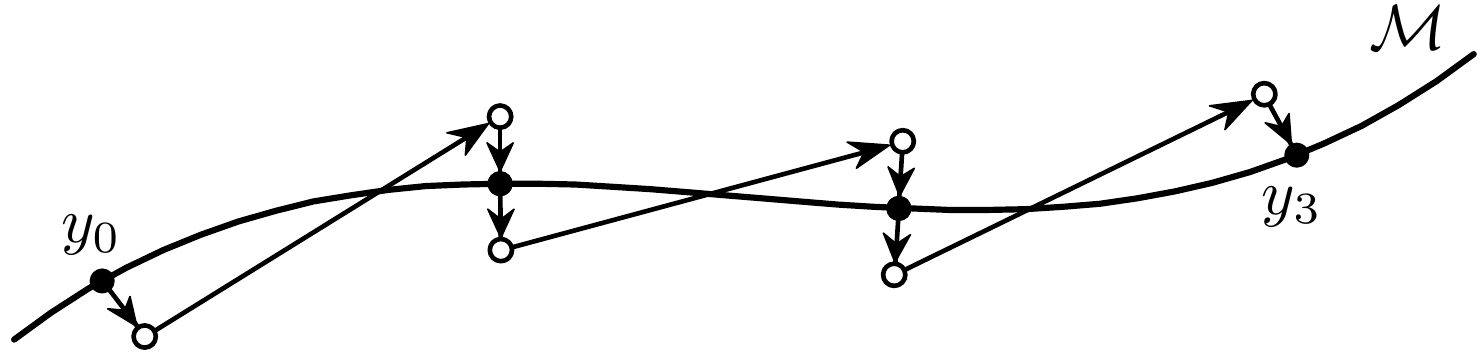}
\caption{The symmetric projection\cite[Fig. 3.3]{Hai:01} }
\label{figb}
\end{center}
\end{figure}

\section{Test problem}

We consider the minimum-time motion over the unit upper hemisphere from a state
$(x_0,y_0,z_0)$ to a state $(x_f,y_f,z_f)$
with inequality constraints.

The system dynamics is governed by the system of differential equations
\begin{equation}\label{e10}
\frac{d}{dt}\left[\begin{array}{c}x\\y\\z\end{array}\right]=
\left[\begin{array}{c}z\cos u\\z\sin u\\-x\cos u-y\sin u\end{array}\right].
\end{equation}
The control variable $u$ is subject to an inequality constraint. Namely,
the control $u$ always stays within the band $c_{u}-r_{u}\leq u\leq c_{u}+r_{u}$.
Following \cite{Oht:04} we introduce a slack variable $u_s$ and replace
the inequality constraint by the equality constraint
\[
C(u,u_d)=(u-c_{u})^2+u_s^2-r_{u}^2=0. 
\]
In order to guarantee that the state passes through the point
$(x_f,y_f,z_f)$ at time $t=t_f$, we impose three terminal constraints
of the form
\[
\psi(x,y,z,p)=\left[\begin{array}{c}x-x_f\\y-y_f\\z-z_f\end{array}\right]=0. 
\]
The objective is to minimize the cost function
\[
J=\phi(p)+\int_{t_0}^{t_f}L(x,y,z,u,u_s,p)dt',
\]
where
\[
\phi(p)=p=t_f-t_0,\quad L(x,y,z,u,u_s,p)=-w_su_s.
\]
The term $\phi(p)$ is responsible for the shortest time to destination,
and the function $L$ serves for stabilization of the slack variable $u_s$.

System (\ref{e10}) possesses the first integral
\[
x^2+y^2+z^2=const. 
\]
Our initial condition lies on the unit sphere $x_0^2+y_0^2+z_0^2=1$.
Therefore, the manifold $\mathcal{M}$ is the unit sphere with the center at the origin.
To simplify the implementation, we choose all parameters such that the
trajectories lie on the upper hemisphere $z\ge0$. The natural local coordinates
in this case are $x$ and $y$, and the local parametrization is given by
$\alpha(x,y)=[x,y,\sqrt{1-x^2-y^2}]^T$.
The system (\ref{e9}) for (\ref{e10}) will be
\begin{equation}\label{e11}
\frac{d}{dt}\left[\begin{array}{c}x\\y\end{array}\right]=
\left[\begin{array}{c}\sqrt{1-x^2-y^2}\cos u\\\sqrt{1-x^2-y^2}\sin u\end{array}\right].
\end{equation}

For this particular example and choice of local coordinates the structure preserving
geometric integration solver consists of the explicit Euler method for the components
$x$ and $y$ and the formula $z=\sqrt{1-x^2-y^2}$ for the component $z$.

For convenience, we change the time variable $t$ within the horizon
by the new time $\tau=(t-t_0)/(t_f-t_0)$, which runs over the interval $[0,1]$.

The corresponding discretized finite-horizon problem on a uniform grid $\tau_i$
in the local coordinates $x$, $y$ comprises the following data structures and computations:
\begin{itemize}
\item $\tau_i=i\Delta\tau$, where $i=0,1,\ldots,N$, and $\Delta\tau=1/N$;
\item the participating variables are the state $\left[\begin{array}{c}
x_i\\y_i\end{array}\right]$, the costate $\left[\begin{array}{c}
\lambda_{1,i}\\\lambda_{2,i}\end{array}\right]$, the control
$\left[\begin{array}{c}u_{i}\\u_{si}\end{array}\right]$, the Lagrange multipliers
$\mu_i$ and $\left[\begin{array}{c}\nu_{1}\\\nu_{2}\end{array}\right]$;
\item the state is governed by the model equation
\[
\left\{\begin{array}{l} x_{i+1}=x_i+\Delta\tau\left(p\sqrt{1-x_i^2-y_i^2}\cos u_{i}\right),\\
\,y_{i+1}=y_i+\Delta\tau\left(p\sqrt{1-x^2-y^2}\sin u_{i}\right),\end{array}\right.
\]
where $i=0,1,\ldots,N-1$;
\item the costate is computed by the backward recursion ($\lambda_{1,N}=\nu_1$,
$\lambda_{2,N}=\nu_2$)
\[
\left\{\begin{array}{l} \lambda_{1,i}=\lambda_{1,i+1}
-\Delta\tau p\frac{x_i}{\sqrt{1-x_i^2-y_i^2}}\\
\hspace{2.5em}{}\cdot(\cos u_i \lambda_{1,i+1}+\sin u_i\lambda_{2,i+1}),\\
\lambda_{2,i} = \lambda_{2,i+1}
-\Delta\tau p\frac{y_i}{\sqrt{1-x_i^2-y_i^2}}\\
\hspace{2.5em}{}\cdot(\cos u_i \lambda_{1,i+1}+\sin u_i\lambda_{2,i+1}),\end{array}\right.
\]
where $i=N-1,N-2,\ldots,0$;
\item the nonlinear equation $F(U,x_0,t_0)=0$, where
\begin{eqnarray*}
\lefteqn{U=[u_0,\ldots,u_{N-1},u_{s,0},\ldots,u_{s,N-1},}\hspace*{8em}\\
&&\mu_0,\ldots,\mu_{N-1},\nu_1,\nu_2,p],
\end{eqnarray*}
has the following rows from the top to bottom:
\[
\left\{\begin{array}{l}       
\Delta\tau p[\sqrt{1-x_i^2-y_i^2}\left(-\sin u_i\lambda_{1,i+1}+
\cos u_i\lambda_{2,i+1}\right)\\
\hspace*{12em}{}+2\left(u_i-c_{u}\right)\mu_i] = 0\end{array}\hspace*{2em}\right.
\]
\[
\left\{\;\;\Delta\tau p\left[2\mu_iu_{si}-w_s\right] = 0\hspace*{12em}\right.
\]
\[
\left\{\;\;\Delta\tau p\left[(u_i-c_{u})^{2}+u_{si}^2-r_{u}^2\right]=0
\right.\hspace*{8em}
\]
\[
\left\{\;\begin{array}{l}x_N-x_f=0\\y_N-y_f=0\end{array}\right.\hspace{17.5em}
\]
\[
\left\{\begin{array}{l}\Delta\tau \{\sum\limits^{N-1}_{i=0}
\sqrt{1-x_i^2-y_i^2}(\cos u_i\lambda_{1,i+1}+\sin u_i\lambda_{2,i+1})\\
\hspace{1em}{}+\mu_i[(u_i-c_u)^2+u_{si}^2-r_u^2]-w_su_{si}\}+1 = 0.\end{array}
\hspace*{4em}\right.
\]
\end{itemize}

\section{Numerical results}

In our numerical experiments, we use the method of the local coordinates
for geometric integration on the unit upper hemisphere.

The number of grid points on the horizon is $N=20$,
the time step of the dynamic system is $\Delta t=0.00625$,
the end points of the computed trajectory are
$(x_0,y_0)=(-0.5,0.5)$ and $(x_f,y_f)=(0.5,0)$.
The constants of the inequality constraint for the control are
$c_u=0.5$ and $r_u=0.1$.
Other parameters are $h=10^{-8}$ and $w_s=0.005$.

We use the GMRES method without restarts. The number of GMRES
iterations does not exceed $k_{\max}=20$, and the
absolute tolerance of the GMRES iterations is $10^{-5}$.

Preconditioning of GMRES accelerates its convergence as
demonstrated on Figures \ref{fig4} and \ref{fig5}.
We used a simple preconditioning strategy as follows. The
exact Jacobian $F_U$ is computed periodically at time instances
with the period $0.2$ seconds. Then the LU factorization of
the Jacobian is used as the preconditioner until the next time
when it is recalculated.

The value of $U$ at $t_0$ is approximated by the Matlab
function \texttt{fsolve} with a special initial guess.

\begin{figure}[ht]
\begin{center}
\includegraphics[width=8.4cm]{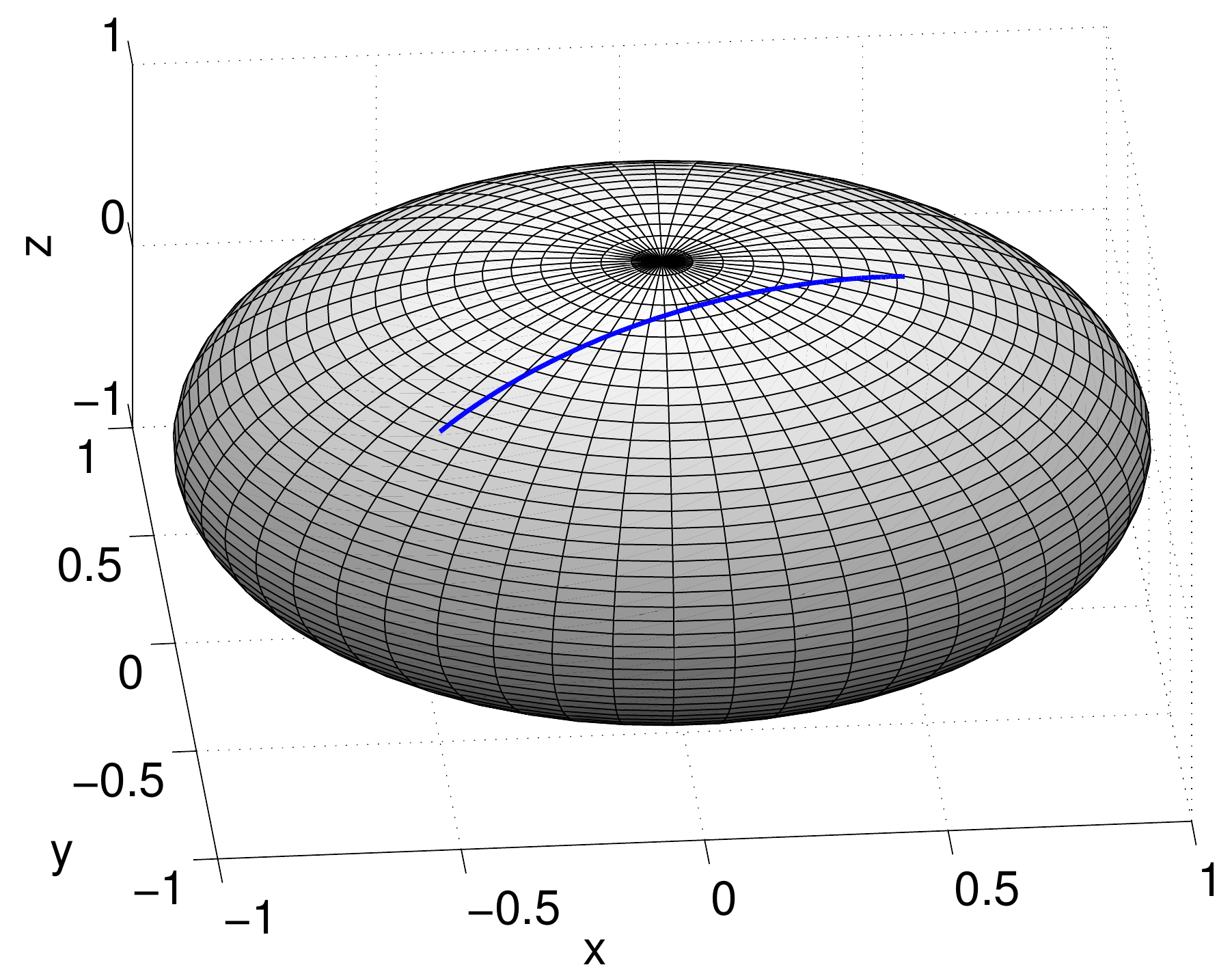}
\caption{3D plot of the trajectory on the unit hemisphere.
Time to destination equals $1.2332$.}
\label{fig1}
\end{center}
\end{figure}

\begin{figure}[ht]
\begin{center}
\includegraphics[width=8.4cm]{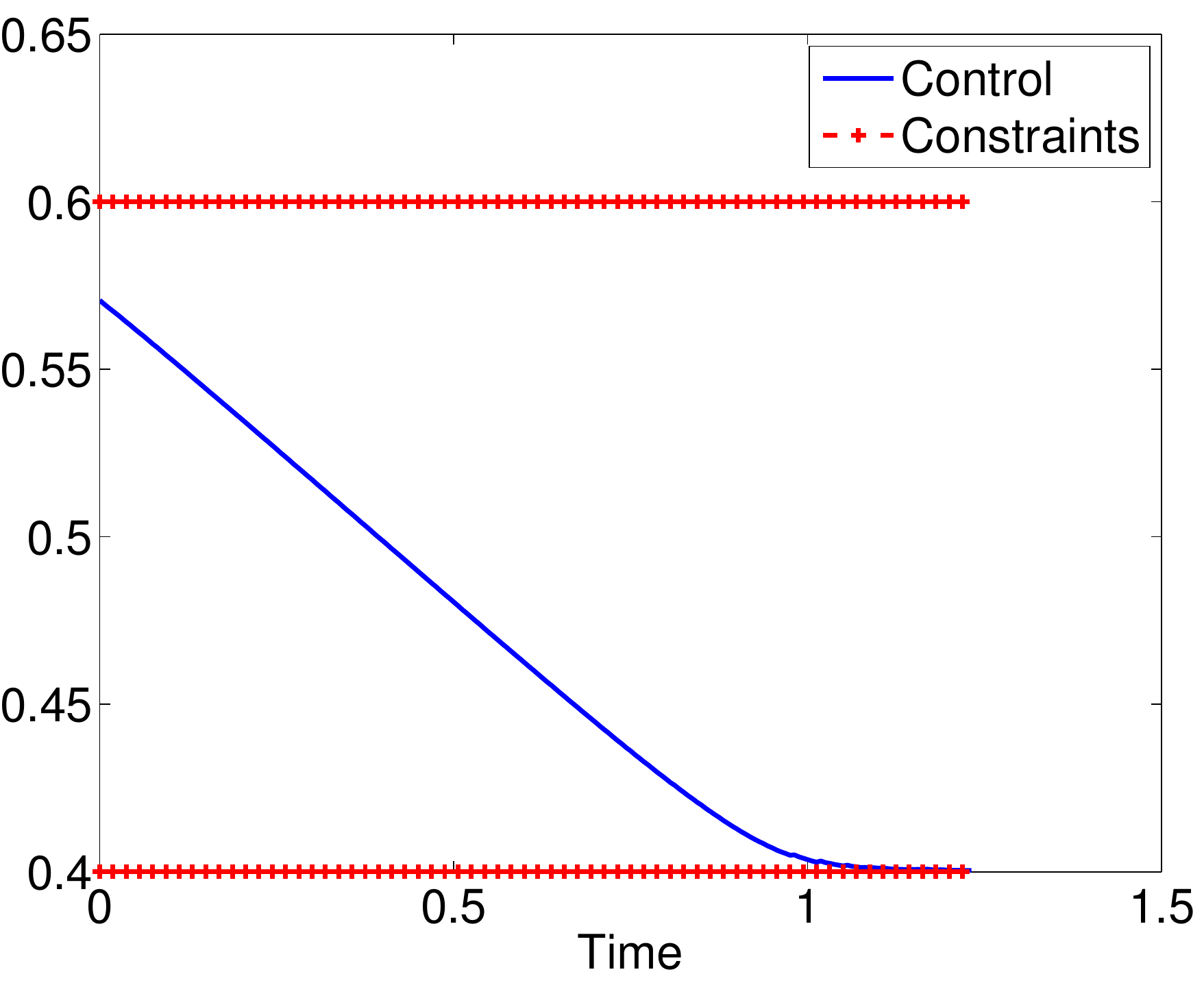}
\caption{The constrained control}
\label{fig2}
\end{center}
\end{figure}

\begin{figure}[ht]
\begin{center}
\includegraphics[width=8.4cm]{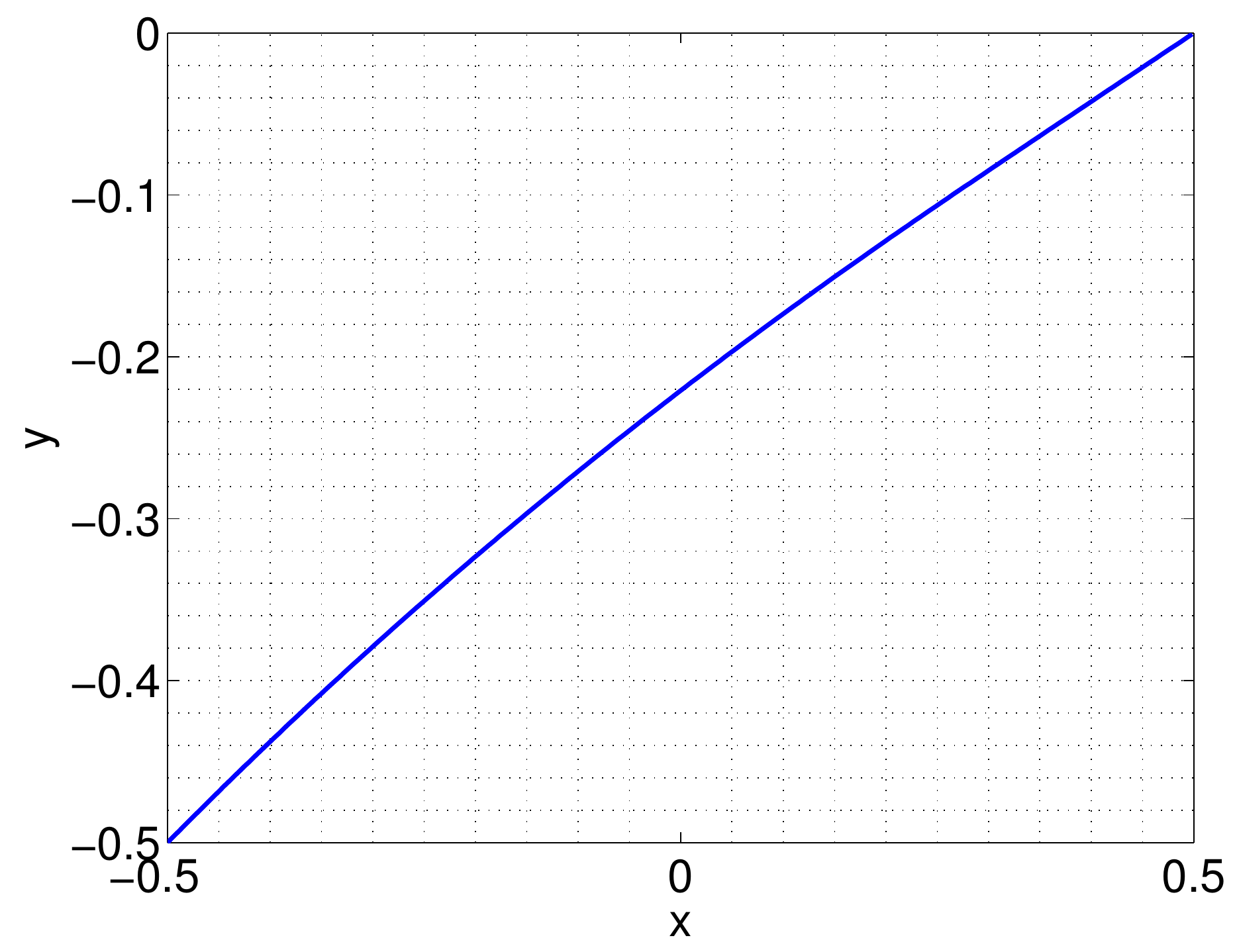}
\caption{2D projection of the trajectory}
\label{fig3}
\end{center}
\end{figure}

\pagebreak\section{Conclusion}

Model predictive control is efficient in dealing with the constraints on
the control and state variables. When the state space of a system lies
on a manifold, it may be profitable to use numerical methods which inherit
this property. For example, numerical algorithms preserving the
symplectic structure of Hamiltonian dynamical systems are much superior
in accuracy compared to the conventional algorithms.
In this note, we show how to incorporate simple modifications
of classical integration methods into our numerical approach to MPC
obtaining an efficient structure-preserving NMPC method. Other
structure-preserving algorithms can be used similarly.

\pagebreak

\begin{figure}[ht]
\begin{center}
\includegraphics[width=7.93cm]{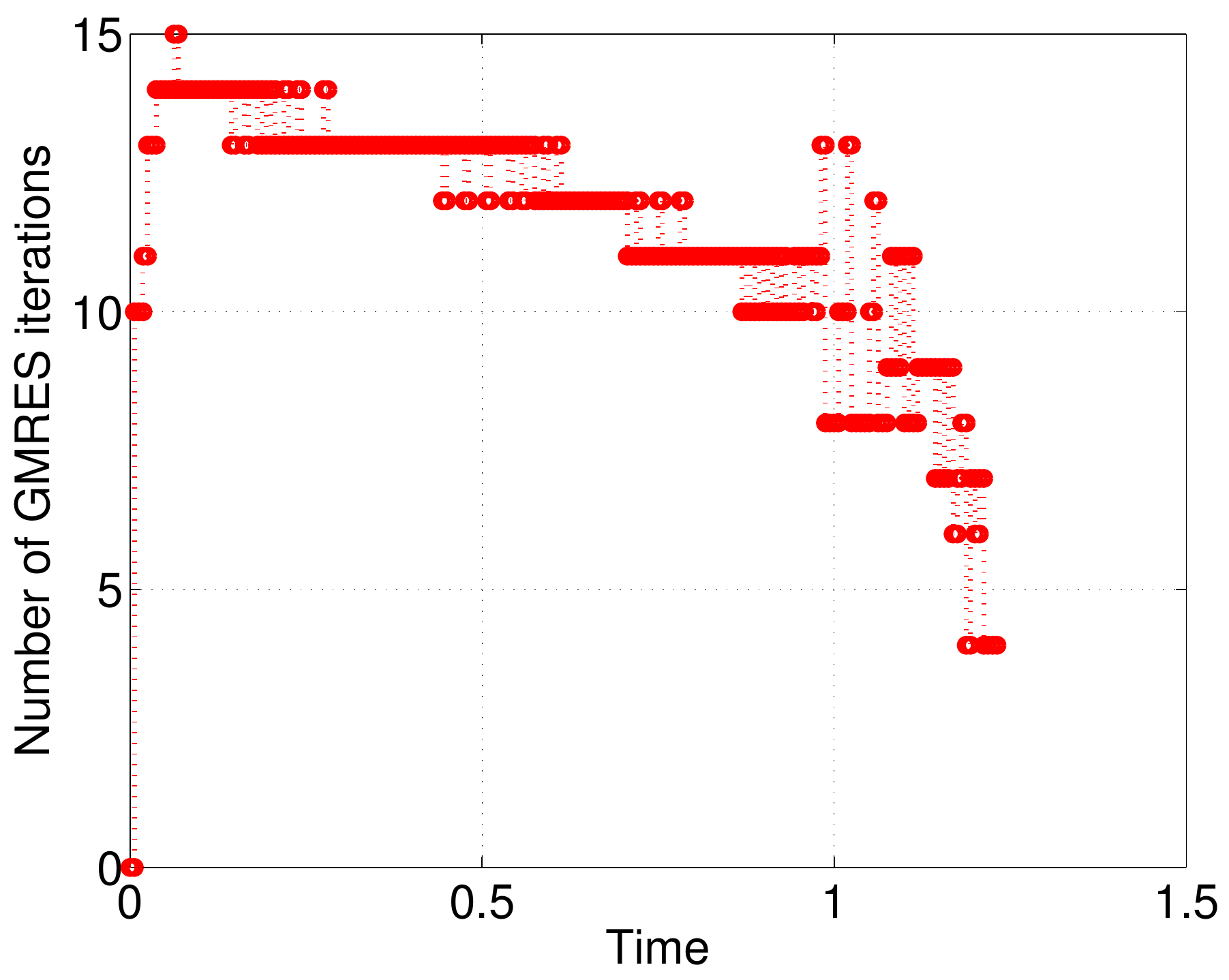}
\caption{Number of GMRES iterations without preconditioning}
\label{fig4}
\end{center}
\end{figure}

\begin{figure}[ht]
\begin{center}
\includegraphics[width=7.93cm]{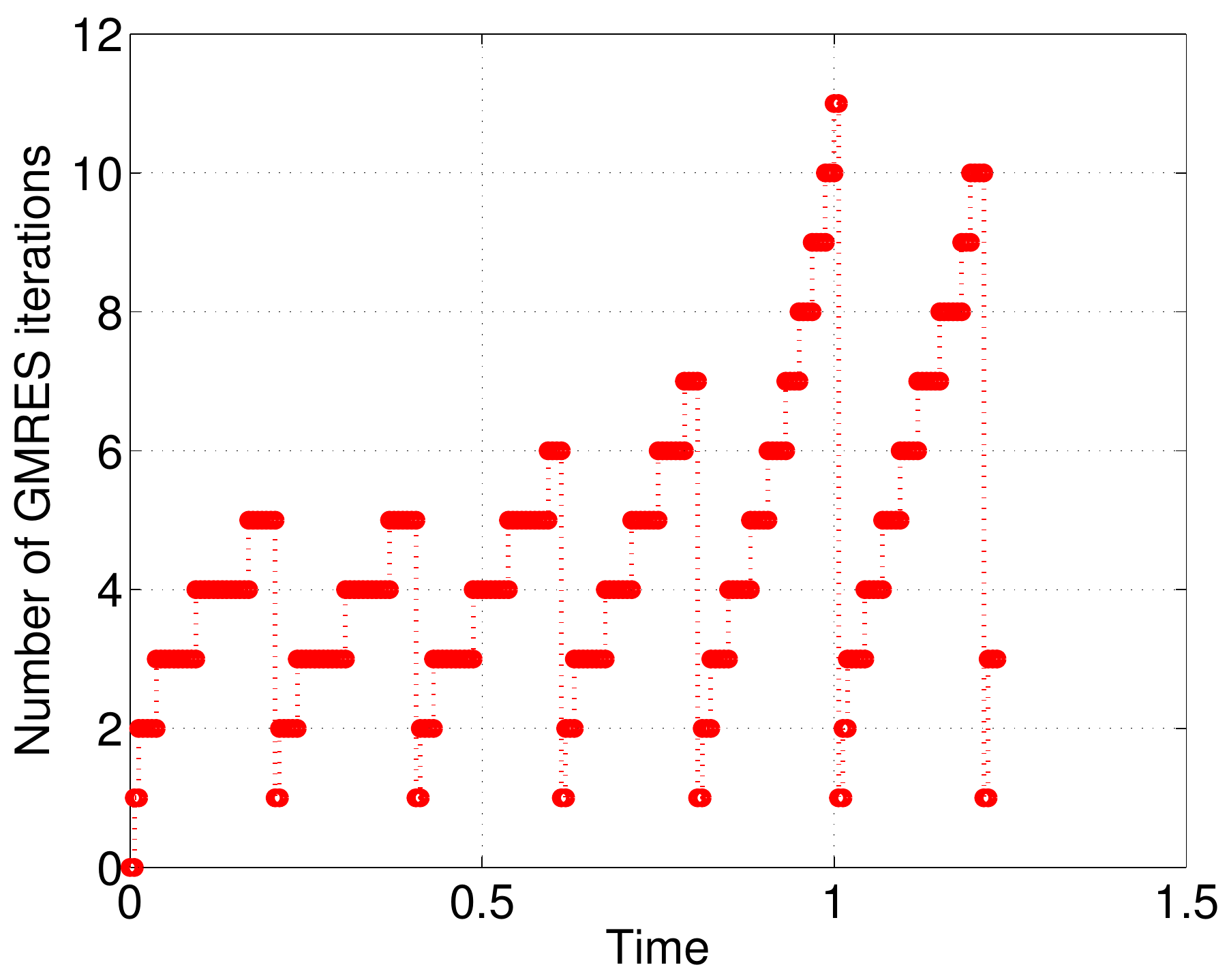}
\caption{Number of GMRES iterations with preconditioning}
\label{fig5}
\end{center}
\end{figure}

\begin{figure}[ht]
\begin{center}
\includegraphics[width=7.93cm]{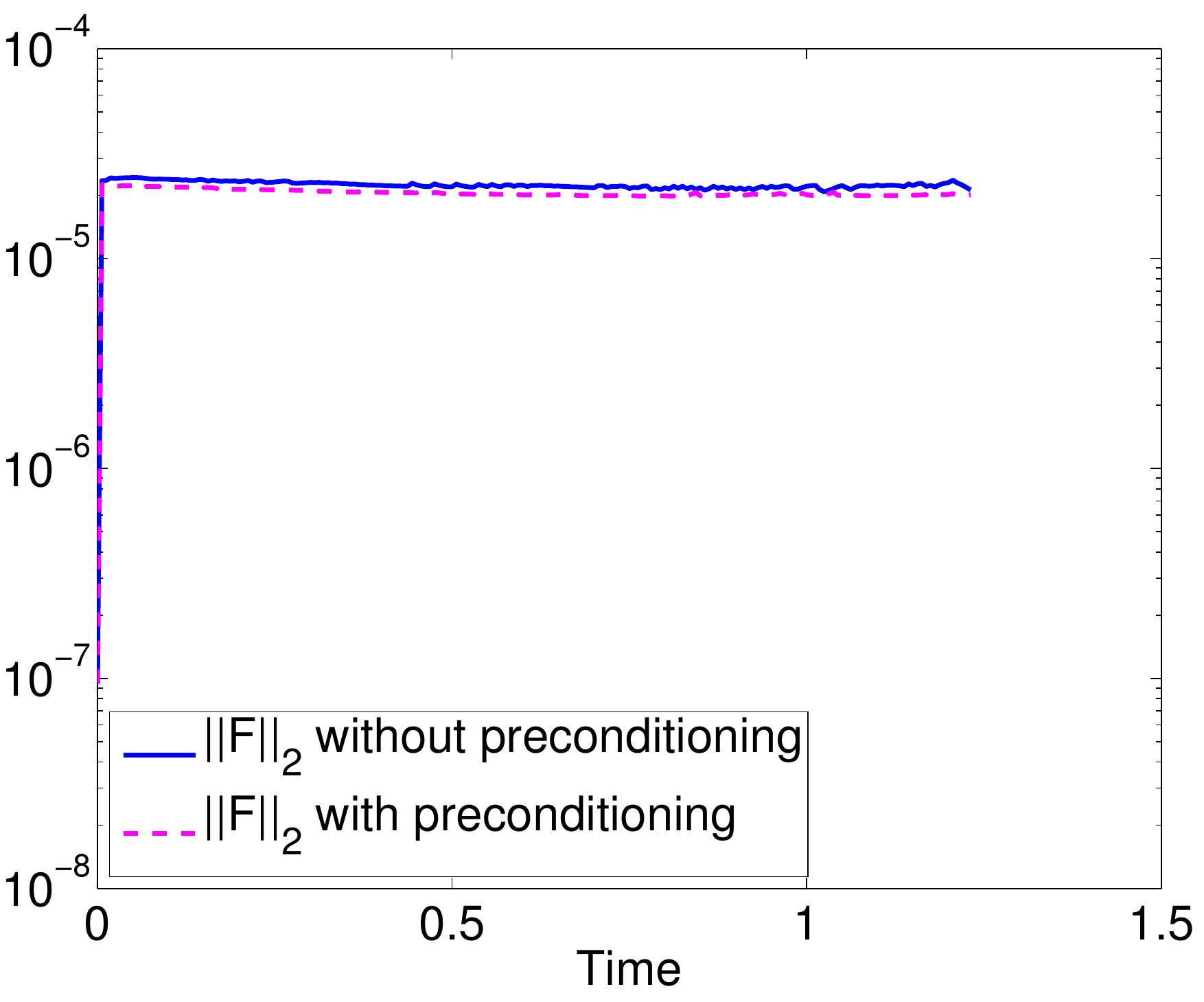}
\caption{Norm of the residual $\|F\|_2$}
\label{fig6}
\end{center}
\end{figure}

\end{document}